# A WEAKNESS IN STRONG LOCALIZATION FOR SINAI'S WALK

By Zhan Shi and Olivier Zindy

*Université Paris VI*

Sinai's walk is a recurrent one-dimensional nearest-neighbor random walk in random environment. It is known for a phenomenon of strong localization, namely, the walk spends almost all time at or near the bottom of deep valleys of the potential. Our main result shows a weakness of this localization phenomenon: with probability one, the zones where the walk stays for the most time can be far away from the sites where the walk spends the most time. In particular, this gives a negative answer to a problem of Erdős and Révész [*Mathematical Structures—Computational Mathematics—Mathematical Modelling* **2** (1984) 152–157], originally formulated for the usual homogeneous random walk.

**1. Introduction.** Let $\omega = (\omega_x, x \in \mathbb{Z})$ be a collection of independent and identically distributed random variables taking values in $(0, 1)$. The distribution of $\omega$ is denoted by $P$. Given the value of $\omega$, we define $(X_n, n \geq 0)$ as a random walk in random environment (RWRE), which is a Markov chain whose distribution is denoted by $P_\omega$. The transition probabilities of $(X_n, n \geq 0)$ are as follows: for $x \in \mathbb{Z}$,

$$P_\omega(X_{n+1} = x+1 | X_n = x) = \omega_x = 1 - P_\omega(X_{n+1} = x-1 | X_n = x).$$

We denote by $\mathbb{P}$ the joint distribution of $(\omega, (X_n))$.

Throughout the paper, we assume that there exists $0 < \delta < \frac{1}{2}$ such that

$$(1.1) \qquad P(\delta \leq \omega_0 \leq 1 - \delta) = 1,$$

and that

$$(1.2) \qquad \mathbb{E}\left[\log\left(\frac{1-\omega_0}{\omega_0}\right)\right] = 0,$$

$$(1.3) \qquad \sigma^2 := \operatorname{Var}\left[\log\left(\frac{1-\omega_0}{\omega_0}\right)\right] > 0.$$

Assumption (1.1) is a commonly adopted technical condition, and can for example be replaced by the existence of exponential moments of $\log(\frac{1-\omega_0}{\omega_0})$. It implies that, $P$-a.s., $|\log(\frac{1-\omega_0}{\omega_0})| \leq M := \log(\frac{1-\delta}{\delta})$. Condition (1.2) ensures, according to Solomon [12], that for $P$-almost all $\omega$, $(X_n)$ is recurrent, that is, it hits any site infinitely often. Finally, (1.3) simply excludes the case of a usual homogeneous random walk.

Recurrent RWRE is known for its slow movement. Indeed, under (1.1)–(1.3), it is proved by Sinai [11] that $X_n/(\log n)^2$ converges in distribution to a nondegenerate limit. Recurrent RWRE will thus be referred to as Sinai's walk. We will from now on assume (1.1)–(1.3).

For an overview of RWRE, see [13]. Although the understanding of one-dimensional RWRE reached a high level in the last decade, there are still some important questions that remain unanswered. See den Hollander [7] for those concerning large deviations.

Let

$$(1.4) \quad \xi(n,x) := \#\{0 \leq i \leq n : X_i = x\}, \qquad n \geq 0, x \in \mathbb{Z},$$

$$(1.5) \quad \mathbb{V}(n) := \left\{x \in \mathbb{Z} : \xi(n,x) = \max_{y \in \mathbb{Z}} \xi(n,y)\right\}, \qquad n \geq 0.$$

In other words, $\xi(n,x)$ records the number of visits at site $x$ by the walk in the first $n$ steps, and $\mathbb{V}(n)$ is the set of sites that are the most visited. Note that $\mathbb{V}(n)$ is not empty. Following Erdős and Révész [4], any element in $\mathbb{V}(n)$ is called a "favorite site."

The basic question we are addressing is: if we know that the walk spends almost all time in $\mathbb{Z}_+$, does it imply that favorite sites would also lie in $\mathbb{Z}_+$?

To formulate the problem more precisely, let us introduce the notion of "positive sequence": a (random) sequence $0 < n_1 < n_2 < \cdots$ of positive numbers is called a "positive sequence" [for the walk $(X_n)$] if

$$(1.6) \quad \lim_{k \to \infty} \frac{\#\{0 \leq i \leq n_k : X_i > 0\}}{n_k} = 1.$$

In other words, the walk spends an overwhelming time in $\mathbb{Z}_+$ along any positive sequence.

PROBLEM 1.1. *Is it true that $\mathbb{P}$-almost surely for any positive sequence $(n_k)$, we have $\mathbb{V}(n_k) \subset \mathbb{Z}_+$ for all large $k$?*

Problem 1.1 was raised by Erdős and Révész [4] (also stated as Problem 10 on page 131 of [8]), originally formulated for the usual homogeneous random walk.

It turns out for the homogeneous walk that the answer is no. Roughly speaking, it is because there is too much "freedom" for the homogeneous walk, so that with probability one, it is possible to find a (random) positive sequence along which the walk does not spend much time in any of the sites of $\mathbb{Z}_+$ – typically, the homogeneous walk makes excursions in $\mathbb{Z}_+$ without spending much time in any sites of $\mathbb{Z}_+$, thus the favorite sites are still in $\mathbb{Z}_-$.

When the environment is random, there is a phenomenon of strong localization [6]; indeed, Sinai's walk spends almost all the time at the bottom of some special zones, called (deep) "valleys." If we know that Sinai's walk spends almost all time in $\mathbb{Z}_+$, then these deep valleys are likely to be located in $\mathbb{Z}_+$, and the favorite sites—which should be located at or near to the bottom of these deep valleys—would also lie in $\mathbb{Z}_+$. In other words, due to strong localization, it looks natural to conjecture that the answer to Problem 1.1 would be yes.

However, things do not go like this. Here is the main result of the paper.

THEOREM 1.2. *Under assumptions* (1.1)–(1.3),

$\mathbb{P}\{\forall \text{ positive sequence } (n_k), \text{ we have } \mathbb{V}(n_k) \subset \mathbb{Z}_+ \text{ for all large } k\} = 0.$

The reason for which the aforementioned heuristics are wrong is that even though Sinai's walk is strongly localized around the bottom of deep valleys, it can happen that a (relatively) big number of sites are around the bottom. In such situations, none of these sites is necessarily favorite, since the visits are shared more or less equally by all these sites.

The main steps in the proof of Theorem 1.2 can be briefly described as follows.

*Step* A. For $P$-almost all environment $\omega$, we define a special sequence, denoted by $(m_k)_{k \geq 1}$. This is the starting point in our construction of a positive sequence $(n_k)$ such that for any $k$, $\mathbb{V}(n_k) \subset \mathbb{Z}_-$.

We mention that the special sequence $(m_k)$ depends only on the environment.

*Step* B. Based on the special sequence $(m_k)$ and on the movement of the walk, we construct in Section 4 another sequence $(n_k)$. We prove that $(n_k)$ is a positive sequence for $(X_n)$, that is, condition (1.6) is satisfied.

*Step* C. Let $(n_k)$ be the positive sequence constructed in Step B. We prove in Section 5 that $\mathbb{P}$-almost surely for all large $k$, $\mathbb{V}(n_k) \subset \mathbb{Z}_-$.

Clearly, Steps B and C together yield Theorem 1.2.

The rest of the paper is organized as follows. In Section 2 we present some elementary facts about Sinai's walk. These facts will be frequently used throughout the paper. A detailed description of Step A is given in Section 3, but the proof of the main result of the section, Proposition 3.1, is postponed to Section 6. Sections 4 and 5 are devoted to Steps B and C, respectively. Finally, in Section 7, we make some comments on the concentration of Sinai's walk.

We use $C_i$ ($1 \leq i \leq 22$) to denote finite and positive constants.

**2. Preliminaries on Sinai's walk.** We list some basic estimates about hitting times and excursions of Sinai's walk.

In the study of Sinai's walk, an important role is played by a process called the potential, denoted by $V = (V(x), x \in \mathbb{Z})$. The potential is a function of the environment $\omega$, and is defined as follows:

$$V(x) := \begin{cases} \sum_{i=1}^{x} \log\left(\frac{1-\omega_i}{\omega_i}\right), & \text{if } x \geq 1, \\ 0, & \text{if } x = 0, \\ -\sum_{i=x+1}^{0} \log\left(\frac{1-\omega_i}{\omega_i}\right), & \text{if } x \leq -1. \end{cases}$$

By (1.1), we have $|V(x) - V(x-1)| \leq M$ for any $x \in \mathbb{Z}$.

2.1. *Hitting times.* For any $x \in \mathbb{Z}$, we define

(2.1) $\quad\quad\quad \tau(x) := \min\{n \geq 1 : X_n = x\}, \quad \min \varnothing := \infty.$

[Attention, if $X_0 = x$, then $\tau(x)$ is the first *return* time to $x$.] Throughout the paper, we write $P_\omega^x(\cdot) := P_\omega(\cdot | X_0 = x)$ (thus $P_\omega^0 = P_\omega$) and write $E_\omega^x$ for expectation with respect to $P_\omega^x$.

It is known ([13], formula (2.1.4)) that for $r < x < s$,

(2.2) $\quad\quad\quad P_\omega^x\{\tau(r) < \tau(s)\} = \sum_{j=x}^{s-1} e^{V(j)} \left(\sum_{j=r}^{s-1} e^{V(j)}\right)^{-1}.$

The next lemma, which gives a simple bound for the expectation of $\tau(r) \wedge \tau(s)$ when the walk starts from a site $x \in (r, s)$, is essentially contained in [6].

LEMMA 2.1. *For any integers $r < s$, we have*

(2.3) $\quad \max_{x \in (r,s) \cap \mathbb{Z}} E_\omega^x[\tau(r) \mathbf{1}_{\{\tau(r) < \tau(s)\}}] \leq (s-r)^2 \exp\left[\max_{r \leq i \leq j \leq s}(V(i) - V(j))\right].$

PROOF. Given $\{\tau(r) < \tau(s)\}$, the walk does not hit site $s$ during time interval $[0, \tau(r)]$. Therefore, $\tau(r)$ under $P_\omega^x\{\cdot | \tau(r) < \tau(s)\}$ is stochastically smaller than the first hitting time of site $r$ by a walk starting from $s$ with a reflecting barrier (to the left) at site $s$. The expected value of this latter random variable is, according to (A1) of [6], bounded by $(s-r)^2 \exp\{\max_{r \leq i \leq j \leq s}(V(i) - V(j))\}$. This yields the lemma. $\square$

We will also use the following estimate borrowed from Lemma 7 of [6]: for $\ell \geq 1$ and $x < y$,

(2.4) $\quad\quad\quad P_\omega^x\{\tau(y) < \ell\} \leq \ell \exp\left(-\max_{x \leq i < y}[V(y-1) - V(i)]\right).$

Looking at the environment backward, we get: for $\ell \geq 1$ and $w < x$,

$$(2.5) \quad P_\omega^x\{\tau(w) < \ell\} \leq \ell \exp\Big(-\max_{w < i \leq x}[V(w+1) - V(i)]\Big).$$

2.2. *Excursions.* We quote some elementary facts about excursions of Sinai's walk (for detailed discussions, see Section 3 of [3]). Let $b \in \mathbb{Z}$ and $x \in \mathbb{Z}$, and consider $\xi(\tau(b), x)$ under $P_\omega^b$. In words, we look at the number of visits to $x$ of the walk (starting from $b$) at the first return to $b$. Then there exist constants $C_1$, $C_2$ and $C_3$ such that

$$(2.6) \quad C_1 e^{-[V(x) - V(b)]} \leq E_\omega^b[\xi(\tau(b), x)] \leq C_2 e^{-[V(x) - V(b)]},$$

and that

$$(2.7) \quad \mathrm{Var}_\omega^b[\xi(\tau(b), x)] \leq C_3 |x - b| \exp\Big(\max_{b \leq y \leq x}[V(y) - V(x)]\Big) e^{-[V(x) - V(b)]},$$

where $\max_{b \leq y \leq x}$ should be replaced by $\max_{x \leq y \leq b}$ if $x < b$.

**3. Step A: a special sequence.** Recall the constant $\delta$ from condition (1.1). We write

$$C_4 := \frac{\delta^3}{2}.$$

For any $j > 0$, we define

$$(3.1) \quad d^+(j) := \min\{n \geq 0 : V(n) \geq j\},$$

$$(3.2) \quad b^+(j) := \min\Big\{n \geq 0 : V(n) = \min_{0 \leq x \leq d^+(j)} V(x)\Big\}.$$

Similarly, we define

$$(3.3) \quad d^-(j) := \max\{n \leq 0 : V(n) \geq j\},$$

$$(3.4) \quad b^-(j) := \max\Big\{n \leq 0 : V(n) = \min_{d^-(j) \leq x \leq 0} V(x)\Big\}.$$

In the next sections, we will be frequently using the following elementary estimates: for any $\varepsilon > 0$, $P$-almost surely for all large $j$,

$$(3.5) \quad j^{2-\varepsilon} \leq |b^\pm(j)| < |d^\pm(j)| \leq j^{2+\varepsilon}.$$

To introduce the announced special sequence in Step A, we define the events [the constant $C_5$ will be defined in (6.8)]:

$$(3.6) \quad E_1^+(j) := \{-2j \leq V(b^+(j)) \leq -j\},$$

$$(3.7) \quad E_2^+(j) := \Big\{\max_{0 \leq x \leq y \leq b^+(j)}[V(y) - V(x)] \leq \frac{j}{4}\Big\},$$

$$(3.8) \quad E_3^+(j) := \left\{ \max_{b^+(j) \leq x \leq y \leq d^+(j)} [V(x) - V(y)] \leq j \right\},$$

$$(3.9) \quad E_4^+(j) := \left\{ \sum_{0 \leq x \leq d^+(j)} e^{-[V(x) - V(b^+(j))]} \geq C_4 \log \log j \right\}$$

and

$$(3.10) \quad E_1^-(j) := \{ V(b^-(j)) \leq -3j \},$$

$$(3.11) \quad E_2^-(j) := \left\{ \max_{b^-(j) \leq x \leq 0} V(x) \geq \frac{j}{3} \right\},$$

$$(3.12) \quad E_3^-(j) := \left\{ \max_{b^-(j) \leq x \leq y \leq 0} [V(x) - V(y)] \leq \frac{j}{2} \right\},$$

$$(3.13) \quad E_4^-(j) := \left\{ \frac{j}{3} \leq \max_{d^-(j) \leq x \leq y \leq b^-(j)} [V(y) - V(x)] \leq j \right\},$$

$$(3.14) \quad E_5^-(j) := \left\{ \sum_{d^-(j) \leq x \leq 0} e^{-[V(x) - V(b^-(j))]} \leq 1 + C_5 \right\}.$$

We set

$$(3.15) \quad E^+(j) := \bigcap_{i=1}^{4} E_i^+(j), \qquad E^-(j) := \bigcap_{i=1}^{5} E_i^-(j).$$

In other words, $E_1^+(j)$, $E_2^+(j)$ and $E_3^+(j)$ require $(V(x), 0 \leq x \leq d^+(j))$ to behave "normally" (i.e., without excessive minimum, nor excessive fluctuations), whereas $E_4^+(j)$ requires the potential to have a "relatively large" number of sites near the minimum.

Similarly, $E_1^-(j)$ and $E_2^-(j)$ require $(V(y), d^-(j) \leq y \leq 0)$ to have no excessive extreme values, $E_3^-(j)$ and $E_4^-(j)$ no excessive fluctuations, $E_5^-(j)$ no excessive concentration around the minimum.

Later, we will see that $P\{E_1^+(j) \cap E_2^+(j) \cap E_3^+(j) \cap E^-(j)\}$ is greater than a positive constant, while $P\{E_4^+(j)\}$ tends to 0 (as $j \to \infty$) "sufficiently slowly."

See Figure 1 for an example of $\omega \in E^+(j) \cap E^-(j)$.

For future use, let us note that for $\omega \in E_3^-(j) \cap E_1^+(j) \cap E_2^+(j) \cap E_3^+(j)$, we have

$$(3.16) \quad \max_{b^-(j) \leq x \leq y \leq d^+(j)} [V(x) - V(y)] \leq \frac{5j}{2}.$$

The proof of the following proposition is postponed until Section 6.

PROPOSITION 3.1. *Under assumptions* (1.1)–(1.3), *for P-almost all environment $\omega$, there exists a random sequence $(m_k)$ such that $\omega \in E^+(m_k) \cap E^-(m_k)$ for all $k \geq 1$.*

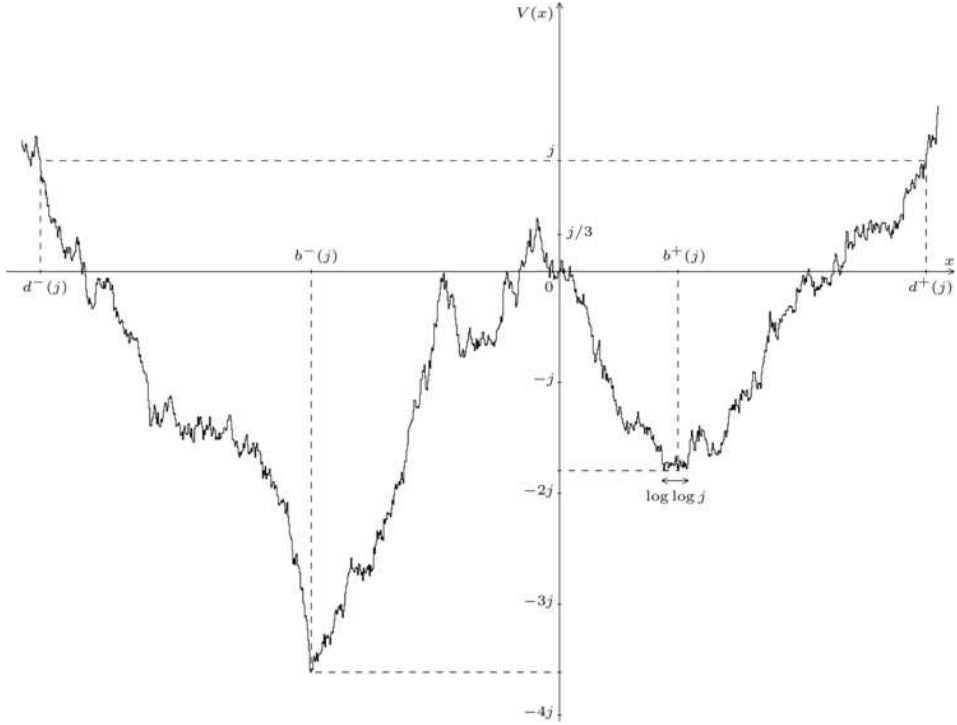

FIG. 1. *Example of $\omega \in E^+(j) \cap E^-(j)$.*

By admitting Proposition 3.1, we will complete Steps B and C in the next two sections.

**4. Step B: a positive sequence.** Let $(m_k)$ be the special sequence introduced in Proposition 3.1. Without loss of generality, we can assume $m_k \geq k^{3k}$ for all $k \geq 1$. For brevity, we write throughout the paper,

(4.1) $\qquad b_k^+ := b^+(m_k), \qquad d_k^+ := d^+(m_k), \qquad \tau_k^+ := \tau(b_k^+),$

(4.2) $\qquad b_k^- := b^-(m_k), \qquad d_k^- := d^-(m_k), \qquad \tau_k^- := \tau(b_k^-).$

We define

(4.3) $$n_k := (1 + (\log k)^{-1/4})\tau_k^-.$$

We prove in this section that $\mathbb{P}$-almost surely, $(n_k)$ is a positive sequence for $(X_n)$, that is, $\frac{1}{n_k}\#\{0 \leq i \leq n_k : X_i > 0\} \to 1$, $\mathbb{P}$-a.s. (as $k \to \infty$).

We start with a few lemmas.

LEMMA 4.1. *We have, P-almost surely, for all large $k$,*

(4.4) $\qquad\qquad\qquad P_\omega\{\tau_k^- < \tau_k^+\} \leq m_k^3 e^{-m_k/12},$

(4.5) $$P_\omega\{\tau(d_k^+) < \tau_k^-\} \leq m_k^3 e^{-m_k/2}.$$

As a consequence, $\mathbb{P}$-almost surely for all large $k$,

(4.6) $$\tau_k^+ < \tau_k^- < \tau(d_k^+).$$

PROOF. By (2.2), $P_\omega\{\tau_k^- < \tau_k^+\} = \sum_{j=0}^{b_k^+-1} e^{V(j)} / \sum_{j=b_k^-}^{b_k^+-1} e^{V(j)}$. Since $\max_{b_k^- \leq j \leq 0} V(j) \geq \frac{m_k}{3}$ [see (3.11)], we have $\sum_{j=b_k^-}^{b_k^+-1} e^{V(j)} \geq e^{m_k/3}$. On the other hand, $\max_{0 \leq j \leq b_k^+} V(j) \leq \frac{m_k}{4}$ according to (3.7). Therefore, $P_\omega\{\tau_k^- < \tau_k^+\} \leq b_k^+ e^{-m_k/12}$. Since $b_k^+ \leq m_k^3$ for large $k$ [see (3.5)], this yields (4.4). The proof of (4.5) is along the same lines, using the fact that $\max_{b_k^- \leq j \leq 0} V(j) \leq \frac{m_k}{2}$ [a consequence of (3.12)].

Since $m_k \geq k$, (4.4) and (4.5) yield, respectively, $\sum_k P_\omega\{\tau_k^- < \tau_k^+\} < \infty$, $\sum_k P_\omega\{\tau(d_k^+) < \tau_k^-\} < \infty$, $P$-almost surely. Now (4.6) follows from the Borel–Cantelli lemma. □

LEMMA 4.2. *Let $A_-(n) := \#\{i : 0 \leq i \leq n, X_i \leq 0\}$. There exists a constant $C_6$ such that $P$-almost surely, for all large $k$,*

(4.7) $$E_\omega[A_-(\tau_k^-)] \leq C_6 (b_k^-)^2 e^{m_k/2}.$$

PROOF. Let $x \in (b_k^-, 0] \cap \mathbb{Z}$. Recall $\xi(n, x)$ from (1.4). Recall that $P_\omega^x(\cdot) := P_\omega(\cdot | X_0 = x)$. Clearly, $P_\omega\{\xi(\tau_k^-, x) = \ell\} = (1 - \pi_x)^{\ell-1} \pi_x$, $\ell \geq 1$, where

$$\pi_x := P_\omega^x\{\tau(x) > \tau_k^-\}$$
(4.8) $$= (1 - \omega_x) P_\omega^{x-1}\{\tau(x) > \tau_k^-\}$$
$$= \frac{1 - \omega_x}{\sum_{j=b_k^-}^{x-1} e^{V(j)-V(x-1)}},$$

the last identity being a consequence of (2.2). In view of assumption (1.1), this yields

$$\frac{1}{\pi_x} \leq C_6 |b_k^-| \exp\left(\max_{b_k^- \leq j \leq i \leq 0} (V(j) - V(i))\right).$$

Since $\max_{b_k^- \leq j \leq i \leq 0} (V(j) - V(i)) \leq \frac{m_k}{2}$, and $E_\omega[\xi(\tau_k^-, x)] = \frac{1}{\pi_x}$, this yields $E_\omega[\xi(\tau_k^-, x)] \leq C_6 |b_k^-| e^{m_k/2}$. Summing over $x \in (b_k^-, 0] \cap \mathbb{Z}$ completes the proof of the lemma. □

REMARK 4.3. *A similar argument shows that for all $x \in [0, b_k^+]$,*

(4.9) $$E_\omega[\xi(\tau_k^+, x)] \leq C_6 b_k^+ e^{m_k/4}.$$

LEMMA 4.4. *For any $k \geq 1$ and $N \geq 1$, we have*

(4.10) $$P_\omega\{\tau_k^+ < \tau_k^- < N\} \leq C_6 e^{-4m_k/3} N.$$

*Furthermore, $P$-almost surely, for all large $k$,*

(4.11) $$E_\omega\left(\frac{1}{\tau_k^-}\mathbf{1}_{\{\tau_k^+ < \tau_k^-\}}\right) \leq C_7 m_k e^{-4m_k/3}.$$

PROOF. We observe that

(4.12) $$P_\omega\{\xi(\tau_k^-, b_k^+) = \ell\} = q_k(1 - \pi_{b_k^+})^{\ell-1}\pi_{b_k^+}, \qquad \ell \geq 1,$$

where $\pi_{b_k^+}$ is as in (4.8), $q_k := P_\omega\{\tau_k^+ < \tau_k^-\}$, and $P_\omega\{\xi(\tau_k^-, b_k^+) = 0\} = 1 - q_k$. Therefore, for any $N \geq 1$, $P_\omega\{1 \leq \xi(\tau_k^-, b_k^+) \leq N\} \leq \pi_{b_k^+} N$. Note that,

(4.13) $$\pi_{b_k^+} \leq \exp\left(V(b_k^+ - 1) - \max_{b_k^- \leq j \leq 0} V(j)\right) \leq C_6 e^{-4m_k/3},$$

the second inequality following from (3.6) and (3.11). In view of the trivial relations $\tau_k^- \geq \xi(\tau_k^-, b_k^+) + 1$ and $\{\tau_k^+ < \tau_k^-\} = \{\xi(\tau_k^-, b_k^+) \geq 1\}$, this implies (4.10).

To prove the second inequality in the lemma, we note that by the strong Markov property, $E_\omega(\frac{\mathbf{1}_{\{\tau_k^+ < \tau_k^-\}}}{1+\xi(\tau_k^-, b_k^+)}) = q_k E_\omega^{b_k^+}(\frac{1}{1+\xi(\tau_k^-, b_k^+)})$. Since $P_\omega^{b_k^+}\{\xi(\tau_k^-, b_k^+) = \ell\} = (1-\pi_{b_k^+})^{\ell-1}\pi_{b_k^+}$, $\ell \geq 1$, it follows that

$$E_\omega\left(\frac{\mathbf{1}_{\{\tau_k^+ < \tau_k^-\}}}{1+\xi(\tau_k^-, b_k^+)}\right) = \frac{q_k \pi_{b_k^+}}{(1-\pi_{b_k^+})^2}\left(\log\left(\frac{1}{\pi_{b_k^+}}\right) - (1 - \pi_{b_k^+})\right)$$

$$\leq \frac{\pi_{b_k^+}}{(1-\pi_{b_k^+})^2}\log\left(\frac{1}{\pi_{b_k^+}}\right).$$

The function $u \mapsto \frac{u}{(1-u)^2}\log(\frac{1}{u})$ is increasing in the (positive) neighborhood of 0. Therefore, by (4.13), $\frac{\pi_{b_k^+}}{(1-\pi_{b_k^+})^2}\log(\frac{1}{\pi_{b_k^+}}) \leq C_7 m_k e^{-4m_k/3}$ (for large $k$). Now (4.11) follows again by means of the trivial inequality $\tau_k^- \geq \xi(\tau_k^-, b_k^+) + 1$. □

LEMMA 4.5. *We have, $\mathbb{P}$-almost surely, for all large $k$,*

(4.14) $$\max_{\tau_k^- \leq i \leq n_k} X_i < 0.$$

PROOF. By the strong Markov property, $\{X(i+\tau_k^-), i\geq 0\}$ is independent (under $P_\omega$) of $\tau_k^-$. Recall that $n_k = (1+(\log k)^{-1/4})\tau_k^-$. We have, for any $\ell \geq 1$,

$$P_\omega\left\{\max_{\tau_k^- \leq i \leq n_k} X_i \geq 0 \Big| \tau_k^- = \ell\right\} = P_\omega^{b_k^-}\left\{\tau(0) \leq \frac{\ell}{(\log k)^{1/4}}\right\} \leq C_6 \frac{\ell}{(\log k)^{1/4}} e^{-3m_k},$$

the last inequality being a consequence of (2.4) [together with (3.10) and (3.11)]. As a consequence,

$$P_\omega\left\{\max_{\tau_k^- \leq i \leq n_k} X_i \geq 0,\, \tau_k^- < \tau(d_k^+)\right\} \leq \frac{C_6}{(\log k)^{1/4}} e^{-3m_k} E_\omega(\tau_k^- \mathbf{1}_{\{\tau_k^- < \tau(d_k^+)\}}).$$

By (2.3) and (3.16),

(4.15)  $$E_\omega(\tau_k^- \mathbf{1}_{\{\tau_k^- < \tau(d_k^+)\}}) \leq (d_k^+ - b_k^-)^2 e^{5m_k/2}.$$

Since $d_k^+ - b_k^- \leq m_k^3$, $P$-almost surely, for all large $k$ [see (3.5)] and since $m_k \geq k$, it follows that

$$\sum_k P_\omega\left\{\max_{\tau_k^- \leq i \leq n_k} X_i \geq 0, \tau_k^- < \tau(d_k^+)\right\} < \infty, \qquad P\text{-a.s.}$$

Recall from (4.6) that $\tau_k^- < \tau(d_k^+)$ $\mathbb{P}$-almost surely, for all large $k$. Lemma 4.5 now follows from the Borel–Cantelli lemma. $\square$

It is now time to complete the argument for Step B by showing that $(n_k)$ is a positive sequence for $(X_n)$.

Combining (4.7) with (4.10) yields that

$$P_\omega\left\{\frac{A_-(\tau_k^-)}{\tau_k^-} \geq e^{-m_k/3}, \tau_k^+ < \tau_k^-\right\}$$
$$\leq P_\omega\{A_-(\tau_k^-) \geq e^{-m_k/3} N\} + P_\omega\{\tau_k^+ < \tau_k^- < N\}$$
$$\leq \frac{C_6(b_k^-)^2 e^{m_k/2}}{e^{-m_k/3} N} + C_6 e^{-4m_k/3} N.$$

Recall that $|b_k^-| \leq m_k^3$ $P$-almost surely, for all large $k$ [see (3.5)]. Choosing $N := e^{m_k}$, and we have, for large $k$,

$$P_\omega\left\{\frac{A_-(\tau_k^-)}{\tau_k^-} \geq e^{-m_k/3}, \tau_k^+ < \tau_k^-\right\} \leq C_8 m_k^6 e^{-m_k/6}.$$

Since $m_k \geq k$, this yields $\sum_k P_\omega\{A_-(\tau_k^-) \geq e^{-m_k/3}\tau_k^-, \tau_k^+ < \tau_k^-\} < \infty$, $P$-almost surely. On the other hand, by (4.6), we have $\tau_k^+ < \tau_k^-$ $\mathbb{P}$-almost surely, for all large $k$. Therefore, the Borel–Cantelli lemma shows that $\mathbb{P}$-almost surely when $k \to \infty$,

$$\frac{A_-(\tau_k^-)}{\tau_k^-} \leq e^{-m_k/3} \to 0. \tag{4.16}$$

Since for large $k$, $A_-(n_k) = A_-(\tau_k^-) + (\log k)^{-1/4} \tau_k^-$ (Lemma 4.5) and $\tau_k^- < n_k$ by definition, we have proved that $\frac{A_-(n_k)}{n_k} \to 0$, $\mathbb{P}$-almost surely. In other words, $(n_k)$ is a positive sequence for the walk.

## 5. Step C: negative favorite sites along a positive sequence.

Let $(n_k)$ be the positive sequence defined in (4.3). In this section, we prove that $\mathbb{P}$-almost surely for all large $k$, $\mathbb{V}(n_k) \subset \mathbb{Z}_-$. As before, we use the notation $b_k^{\pm}$, $d_k^{\pm}$ and $\tau_k^{\pm}$ as in (4.1)–(4.2).

We will prove that $\mathbb{P}$-almost surely, for all large $k$,

$$\xi(n_k, b_k^-) \geq \frac{\tau_k^-}{(\log k)^{1/3}}, \tag{5.1}$$

$$\max_{x \in [1, d_k^+]} \xi(\tau_k^-, x) \leq \frac{\tau_k^-}{(\log k)^{1/2}}. \tag{5.2}$$

Observe that $\mathbb{P}$-almost surely, for all large $k$, we have $\max_{x \in [1, d_k^+]} \xi(\tau_k^-, x) = \max_{x \geq 1} \xi(\tau_k^-, x)$ [by (4.6)], and $\max_{x \geq 1} \xi(\tau_k^-, x) = \max_{x \geq 1} \xi(n_k, x)$ (Lemma 4.5). It is now clear that (5.1) and (5.2) together will complete Step C, and thus the proof of Theorem 1.2.

The rest of the section is devoted to the proof of (5.1) and (5.2). For the sake of clarity, they are proved in distinct subsections.

### 5.1. Proof of (5.1).

Let $T_0^- := \tau_k^-$ and

$$T_j^- = T_j^-(k) := \min\{n > T_{j-1}^- : X_n = b_k^-\}, \qquad j = 1, 2, \ldots.$$

We define, for any $j \geq 1$,

$$Y_j^-(x) := \xi(T_j^-, x) - \xi(T_{j-1}^-, x), \qquad x \in \mathbb{Z},$$

$$Z_j^- := \sum_{x \in (d_k^-, 0]} Y_j^-(x).$$

By the strong Markov property, $(Z_j^-, j \geq 1)$ is a sequence of i.i.d. random variables (under $P_\omega$). Recall that $n_k = (1 + (\log k)^{-1/4}) \tau_k^-$. Let $\ell \geq 1$. By the strong Markov property, $P_\omega\{\xi(n_k, b_k^-) < \frac{\ell}{(\log k)^{1/3}} | \tau_k^- = \ell\} = P_\omega^{b_k^-}\{\xi(\frac{\ell}{(\log k)^{1/4}}, b_k^-) < \frac{\ell}{(\log k)^{1/3}}\}$. Under probability $P_\omega^{b_k^-}$, if $\tau(d_k^-) \wedge \tau(0) > \frac{\ell}{(\log k)^{1/4}}$, then the walk stays in $(d_k^-, 0)$ during time interval $[0, \frac{\ell}{(\log k)^{1/4}}]$; if moreover $\xi(\frac{\ell}{(\log k)^{1/4}}, b_k^-) < \frac{\ell}{(\log k)^{1/3}}$, then $\sum_{j=1}^{\ell/(\log k)^{1/3}} Z_j^- \geq \frac{\ell}{(\log k)^{1/4}}$. Accordingly,

$$P_\omega\Big\{\xi(n_k, b_k^-) < \frac{\ell}{(\log k)^{1/3}} \Big| \tau_k^- = \ell\Big\}$$

$$\leq P_\omega\Big\{\sum_{j=1}^{\ell/(\log k)^{1/3}} Z_j^- \geq \frac{\ell}{(\log k)^{1/4}}\Big\} + P_\omega^{b_k^-}\Big\{\tau(d_k^-) \wedge \tau(0) \leq \frac{\ell}{(\log k)^{1/4}}\Big\}.$$

By (2.6), $E_\omega(Z_j^-) \leq C_2 \sum_{x \in (d_k^-, 0]} e^{-[V(x) - V(b_k^-)]}$, which, according to (3.14), is bounded by $C_2(1 + C_5) =: C_9$. Therefore,

$$P_\omega\Big\{\sum_{j=1}^{\ell/(\log k)^{1/3}} Z_j^- \geq \frac{\ell}{(\log k)^{1/4}}\Big\}$$

$$\leq P_\omega\Big\{\sum_{j=1}^{\ell/(\log k)^{1/3}} (Z_j^- - E_\omega Z_j^-) \geq ((\log k)^{-1/4} - C_9(\log k)^{-1/3})\ell\Big\}$$

$$\leq \frac{\text{Var}_\omega(Z_1^-)}{((\log k)^{-1/4} - C_9(\log k)^{-1/3})^2 (\log k)^{1/3} \ell}.$$

We have $\text{Var}_\omega(Z_1^-) \leq |d_k^-| \sum_{x \in (d_k^-, 0]} \text{Var}_\omega(Y_1^-(x))$. By (2.7) and (3.12)–(3.13), $\text{Var}_\omega(Y_1^-(x))$ is bounded by $C_3|d_k^-|e^{m_k}$ for all $x \in (d_k^-, 0]$. Thus $\text{Var}_\omega(Z_1^-) \leq C_3(d_k^-)^2 e^{m_k}$. Accordingly, for large $k$,

$$P_\omega\Big\{\sum_{j=1}^{\ell/(\log k)^{1/3}} Z_j^- \geq \frac{\ell}{(\log k)^{1/4}}\Big\} \leq C_{10} \frac{(\log k)^{1/6}(d_k^-)^2 e^{m_k}}{\ell}.$$

We now estimate $P_\omega^{b_k^-}\{\tau(d_k^-) \wedge \tau(0) \leq \frac{\ell}{(\log k)^{1/4}}\}$. There is nothing to estimate if $\ell < (\log k)^{1/4}$, so let us assume $\ell \geq (\log k)^{1/4}$. By (2.5) and (3.10),

$$P_\omega^{b_k^-}\Big\{\tau(d_k^-) \leq \frac{\ell}{(\log k)^{1/4}}\Big\} \leq \Big(\frac{\ell}{(\log k)^{1/4}} + 1\Big) e^{-[V(d_k^- + 1) - V(b_k^-)]}$$

$$\leq \frac{C_6 \ell}{(\log k)^{1/4}} e^{-4 m_k},$$

whereas by (2.4) and (3.10),

$$P_\omega^{b_k^-}\Big\{\tau(0) \leq \frac{\ell}{(\log k)^{1/4}}\Big\} \leq \Big(\frac{\ell}{(\log k)^{1/4}} + 1\Big) e^{-[V(-1) - V(b_k^-)]}$$

$$\leq \frac{C_6 \ell}{(\log k)^{1/4}} e^{-3 m_k}.$$

Thus, for all $\ell \geq 1$,

$$P_\omega^{b_k^-}\left\{\tau(d_k^-) \wedge \tau(0) \leq \frac{\ell}{(\log k)^{1/4}}\right\} \leq \frac{2C_6 \ell}{(\log k)^{1/4}} e^{-3m_k}$$

$$=: \frac{C_{11}\ell}{(\log k)^{1/4}} e^{-3m_k}.$$

As a consequence, we have proved that

$$P_\omega\left\{\xi(n_k, b_k^-) \leq \frac{\ell}{(\log k)^{1/3}}\Big|\tau_k^- = \ell\right\}$$

$$\leq \frac{C_{10}(\log k)^{1/6}(d_k^-)^2 e^{m_k}}{\ell} + \frac{C_{11}\ell}{(\log k)^{1/4}} e^{-3m_k}.$$

Therefore,

$$P_\omega\{\xi(n_k, b_k^-) \leq (\log k)^{-1/3}\tau_k^-, \tau_k^+ < \tau_k^- < \tau(d_k^+)\}$$

$$\leq C_{10}(\log k)^{1/6}(d_k^-)^2 e^{m_k} E_\omega\left(\frac{\mathbf{1}_{\{\tau_k^+ < \tau_k^-\}}}{\tau_k^-}\right)$$

$$+ \frac{C_{11}}{(\log k)^{1/4}} e^{-3m_k} E_\omega(\tau_k^- \mathbf{1}_{\{\tau_k^- < \tau(d_k^+)\}}).$$

The two expectations, $E_\omega(\frac{1}{\tau_k^-}\mathbf{1}_{\{\tau_k^+ < \tau_k^-\}})$ and $E_\omega(\tau_k^-\mathbf{1}_{\{\tau_k^- < \tau(d_k^+)\}})$, are estimated by means of (4.11) and (4.15), respectively. We have therefore proved that, for large $k$, $P_\omega\{\xi(n_k, b_k^-) \leq (\log k)^{-1/3}\tau_k^-, \tau_k^+ < \tau_k^- < \tau(d_k^+)\}$ is bounded by

$$C_{10}C_7(\log k)^{1/6}(d_k^-)^2 m_k e^{-m_k/3} + \frac{C_{11}}{(\log k)^{1/4}}(d_k^+ - b_k^-)^2 e^{-m_k/2}.$$

Since $|d_k^-| \leq m_k^3$ and $d_k^+ - b_k^- \leq m_k^3$, $P$-almost surely, for all large $k$ [see (3.5)], and since $m_k \geq k$, this implies

$$\sum_k P_\omega\{\xi(n_k, b_k^-) \leq (\log k)^{-1/3}\tau_k^-, \tau_k^+ < \tau_k^- < \tau(d_k^+)\} < \infty, \qquad P\text{-a.s.}$$

The proof of (5.1) is now completed by means of the Borel–Cantelli lemma and (4.6).

5.2. *Proof of* (5.2). The proof of (5.2) bears many similarities to the proof of (5.1), the basic idea being again via excursions.

Let $T_0^+ := \tau_k^+$ and

$$T_j^+ = T_j^+(k) := \inf\{n > T_{j-1}^+ : X_n = b_k^+\}, \qquad j = 1, 2, \ldots.$$

We write, for any $j \geq 1$,

$$Y_j^+(y) := \xi(T_j^+, y) - \xi(T_{j-1}^+, y), \qquad y \in \mathbb{Z},$$

$$Z_j^+ := \sum_{y \in [1, d_k^+]} Y_j^+(y).$$

Let $M = M(k) := \max\{j : T_j^+ < \tau_k^-\}$. In other words, $M$ denotes the number of excursions (away from $b_k^+$) completed by the walk before hitting $b_k^-$.

Let $x \in [1, d_k^+]$. We have $\xi(\tau_k^-, x) \leq \xi(\tau_k^+, x) + \sum_{j=1}^{M+1} Y_j^+(x)$ and $\#\{i \leq \tau_k^- : X_i \geq 0\} \geq \sum_{j=1}^{M} Z_j^+$. Note that $\{M \geq 1\} = \{\tau_k^+ < \tau_k^-\}$. Therefore, for any $\ell \geq 1$ and $k_r := \ell 2^r$,

$$p(x) := P_\omega\{(\log k)^{1/2} \xi(\tau_k^-, x) > \#\{i \leq \tau_k^- : X_i \geq 0\}, \tau_k^+ < \tau_k^-\}$$

$$\leq P_\omega\{1 \leq M < \ell\} + \sum_{r=0}^{\infty} P_\omega\left\{(\log k)^{1/2} \xi(\tau_k^-, x) > \sum_{j=1}^{M} Z_j^+, k_r \leq M < k_{r+1}\right\}$$

$$\leq P_\omega\{1 \leq M < \ell\} + P_\omega\{\xi(\tau_k^+, x) > \ell\} + \sum_{r=0}^{\infty} I_r,$$

where

$$I_r := P_\omega\left\{(\log k)^{1/2} \ell + (\log k)^{1/2} \sum_{j=1}^{k_{r+1}} Y_j^+(x) > \sum_{j=1}^{k_r} Z_j^+\right\}.$$

By (4.12), we have $P_\omega\{1 < \xi(\tau_k^-, b_k^+) \leq \ell\} \leq \pi_{b_k^+} \ell$, whereas by (4.9), $P_\omega\{\xi(\tau_k^+, x) > \ell\} \leq \frac{1}{\ell} E_\omega[\xi(\tau_k^+, x)] \leq \frac{C_6 b_k^+}{\ell} e^{m_k/4}$. Thus,

(5.3) $$p(x) \leq \pi_{b_k^+} \ell + \frac{C_6 b_k^+}{\ell} e^{m_k/4} + \sum_{r=0}^{\infty} I_r.$$

We now estimate $I_r$. Recall that $Y_j^+(x)$ is the number of visits at site $x$ by an excursion (away from $b_k^+$). According to (2.6), $E_\omega[Y_1^+(x)] \leq C_2 e^{-[V(x)-V(b_k^+)]} \leq C_2$. On the other hand, it follows from (2.6) and then (3.9) that $E_\omega(Z_1^+) \geq C_1 \sum_{y \in [1, d_k^+]} e^{-[V(y)-V(b_k^+)]} \geq C_1 C_4 \log \log m_k$. Since $(\log k)^{1/2} \ell - C_1 C_4 k_r \times \log \log m_k + C_2 (\log k)^{1/2} k_{r+1} \leq -\frac{C_1 C_4}{2} k_r \log \log m_k$ (for large $k$; recalling that $m_k \geq k^{3k}$), we see that, $P$-almost surely, for all large $k$, the probability $I_r$ is bounded (uniformly in all $r \geq 0$) by

$$P_\omega\left\{\sum_{j=1}^{k_r} [Z_j^+ - E_\omega(Z_j^+)] - (\log k)^{1/2} \sum_{j=1}^{k_{r+1}} [Y_j^+(x) - E_\omega(Y_j^+(x))]\right.$$

$$< -\frac{C_1 C_4}{2}(\log \log m_k)k_r\bigg\}$$

$$\leq \frac{8}{(C_1 C_4 \log \log m_k)^2 k_r}[\operatorname{Var}_\omega(Z_1^+) + 2(\log k)\operatorname{Var}_\omega(Y_1^+(x))].$$

By means of (2.7) and (3.7)–(3.8), $\operatorname{Var}_\omega(Y_1^+(x)) \leq C_3 d_k^+ e^{m_k}$; it follows that $\operatorname{Var}_\omega(Z_1^+) \leq d_k^+ \sum_{y \in [1, d_k^+]} \operatorname{Var}_\omega(Y_1^+(y)) \leq C_3 (d_k^+)^3 e^{m_k}$. Accordingly,

$$I_r \leq \frac{8 C_3 d_k^+ [(d_k^+)^2 + 2\log k] e^{m_k}}{(C_1 C_4 \log \log m_k)^2 k_r}.$$

Plugging this into (5.3), and using the fact that $\sum_r k_r^{-1} = 2\ell^{-1}$, we get that for any $\ell \geq 1$,

$$\max_{x \in [1, d_k^+]} p(x) \leq \pi_{b_k^+} \ell + \frac{C_6 b_k^+}{\ell} e^{m_k/4} + C_{12}\frac{d_k^+[(d_k^+)^2 + 2\log k]e^{m_k}}{(\log \log m_k)^2 \ell}.$$

Recall from (4.13) that $\pi_{b_k^+} \leq C_6 e^{-4m_k/3}$. Now we choose $\ell := e^{5m_k/4}$, to see that by (3.5),

$$\sum_k d_k^+ \max_{x \in [1, d_k^+]} p(x) < \infty, \qquad P\text{-a.s.}$$

This implies that $\sum_k P_\omega\{(\log k)^{1/2} \max_{x \in [1, d_k^+]} \xi(\tau_k^-, x) > \#\{i \leq \tau_k^- : X_i \geq 0\}, \tau_k^+ < \tau_k^-\} < \infty$, $P$-almost surely. This implies (5.2) by an application of the Borel–Cantelli lemma and (4.6).

**6. Proof of Proposition 3.1.** We now prove that, for $P$-almost all environment $\omega$, there exists a sequence $(m_k)$ such that $\omega \in E^+(m_k) \cap E^-(m_k)$, $\forall k \geq 1$, where $E^+(m_k)$ and $E^-(m_k)$ are defined in (3.15).

Let $j_k := k^{3k}$ $(k \geq 1)$ and $\mathcal{F}_{j_{k-1}} := \sigma\{V(x), 0 \leq x \leq d^+(j_{k-1})\}$.

Recall that $(E_j^+)$ and $(E_j^-)$ are independent events. If we are able to show that

(6.1) $$\sum_k P\{E^+(j_k)|\mathcal{F}_{j_{k-1}}\} = \infty, \qquad P\text{-a.s.},$$

and that for some $C_- > 0$ and all large $j$,

(6.2) $$P\{E^-(j)\} \geq C_-,$$

then Lévy's Borel–Cantelli lemma ([10], page 518) will tell us that with positive probability, there are infinitely many $k$ such that $\omega \in E^+(j_k) \cap E^-(j_k)$. An application of the Hewitt–Savage zero–one law ([5], Theorem IV.6.3) will then yield Proposition 3.1.

The rest of the section is devoted to the proof of (6.1) and (6.2), presented in distinct subsections.

6.1. *Proof of* (6.1). Recall that $|V(x) - V(x-1)| \leq M = \log\frac{1-\delta}{\delta}$ for any $x \in \mathbb{Z}$.

To bound $P\{E^+(j_k)|\mathcal{F}_{j_{k-1}}\}$ from below, we start with the trivial inequality $E^+(j_k) \supset E^+(j_k) \cap B^+(j_{k-1})$, for any set $B^+(j_{k-1})$. We choose

$$B^+(j_{k-1}) := \left\{\inf_{0 \leq y \leq d^+(j_{k-1})} V(y) \geq -j_{k-1}\log^2 j_{k-1}\right\}.$$

Clearly, $B^+(j_{k-1})$ is $\mathcal{F}_{j_{k-1}}$-measurable. Moreover, on $B^+(j_{k-1}) \cap E^+(j_k)$, we have $d^+(j_{k-1}) \leq b^+(j_k)$.

Recall that $E^+(j_k) = \bigcap_{i=1}^4 E_i^+(j_k)$. Let

$$F_2^+(j_k) := \left\{\max_{0 \leq x \leq y \leq b^+(j_k)}[V(y) - V(x)] \leq \frac{j_k}{4} - j_{k-1}\log^2 j_{k-1} - j_{k-1} - M\right\}.$$

We consider

$$F^+(j_k) := E_1^+(j_k) \cap E_3^+(j_k) \cap E_4^+(j_k) \cap F_2^+(j_k).$$

Since $V(d^+(j_{k-1})) \in I_{j_{k-1}} := [j_{k-1}, j_{k-1} + M]$, we have, by applying the strong Markov property at $d^+(j_{k-1})$,

$$P\{E^+(j_k)|\mathcal{F}_{j_{k-1}}\} \geq \left(\inf_{z \in I_{j_{k-1}}} P_z\{F^+(j_k)\}\right)\mathbf{1}_{B^+(j_{k-1})},$$

where $P_z(\cdot) := P(\cdot|V(0) = z)$, for any $z \in \mathbb{R}$; thus $P = P_0$. (Strictly speaking, we should be working in a canonical space for $V$, with $P_z$ defined as the image measure of $P$ under translation.)

Clearly, $\mathbf{1}_{B^+(j_{k-1})} = 1$, $P$-almost surely for all large $k$. Therefore, the proof of (6.1) boils down to showing the existence of a positive constant $C^+$ such that $P$-almost surely for all large $k$,

(6.3) $$\inf_{z \in I_{j_{k-1}}} P_z\{F^+(j_k)\} \geq \frac{C^+}{k}.$$

Let, for any Borel set $A \subset \mathbb{R}$,

$$d^+(A) := \inf\{i \geq 0 : V(i) \in A\}.$$

A simple martingale argument yields that, whenever $x < y < z$,

(6.4) $$P_y\{d^+([z,\infty)) < d^+((-\infty,x])\} \geq \frac{y-x}{z-x+M},$$

(6.5) $$P_y\{d^+((-\infty,x)) < d^+([z,\infty))\} \geq \frac{z-y}{z-x+M}.$$

We now proceed to prove (6.3). Let

$$a_\ell := -2j_k + 3M\ell, \qquad G_1^+(j_k,\ell) := \{a_\ell \leq V(b^+(j_k)) < a_{\ell+1}\}.$$

Then

$$P_z\{F^+(j_k)\} = P_z\{E_1^+(j_k), F_2^+(j_k), E_3^+(j_k), E_4^+(j_k)\}$$

(6.6)
$$\geq \sum_{\ell=0}^{\lfloor j_k/(3M) \rfloor - 1} P_z\{G_1^+(j_k, \ell), F_2^+(j_k), E_3^+(j_k), E_4^+(j_k)\}$$

$$=: \sum_{\ell=0}^{\lfloor j_k/(3M) \rfloor - 1} P_{k,\ell}^+.$$

Let $L(k,\ell) := \#\{0 \leq i \leq d^+(j_k) : V(i) \in [a_\ell, a_{\ell+1})\}$. On $G_1^+(j_k, \ell)$, we clearly have $e^{-3M} L(k,\ell) \leq \sum_{0 \leq x \leq d^+(j_k)} e^{-[V(x)-V(b^+(j_k))]}$. Thus

$$P_{k,\ell}^+ \geq P_z\{G_1^+(j_k, \ell), F_2^+(j_k), E_3^+(j_k), e^{-3M} L(k,\ell) \geq C_4 \log \log j_k\}$$

$$\geq P_z\{G_1^+(j_k, \ell), F_2^+(j_k), E_3^+(j_k), L(k,\ell) \geq \tfrac{1}{2} \log \log j_k\},$$

the last inequality following from the values of

$$M := \log \frac{1-\delta}{\delta} \quad \text{and} \quad C_4 := \frac{\delta^3}{2}.$$

We define $T_0 := 0$, and by induction,

$$\tau_p := \min\{i \geq T_{p-1} : V(i) < a_{\ell+1}\},$$
$$T_p := \min\{i \geq \tau_p : V(i) \geq a_{\ell+1}\}, \qquad p = 1, 2, \ldots.$$

Let

$$\alpha = \alpha(k) := \lfloor \tfrac{1}{2} \log \log j_k \rfloor,$$
$$\widetilde{T} := \min\{i \geq \tau_1 : V(i) \geq a_{\ell+2}\}.$$

Since $G_1^+(j_k, \ell) \cap \{L(k,\ell) \geq \alpha\} \supset \{\tau_\alpha < \widetilde{T} < d^+(j_k) < d^+((-\infty, a_\ell])\}$, we have

$$P_{k,\ell}^+ \geq P_z\{\tau_\alpha < \widetilde{T} < d^+(j_k) < d^+((-\infty, a_\ell]), F_2^+(j_k), E_3^+(j_k)\}.$$

Consider now the events

$$H_2^+(j_k) := \left\{ \max_{0 \leq x \leq y \leq \tau_1} [V(y) - V(x)] \leq \frac{j_k}{5} \right\},$$

$$H_3^+(j_k) := \left\{ \max_{\widetilde{T} \leq x \leq y \leq d^+(j_k)} [V(x) - V(y)] \leq j_k \right\}.$$

We have, for large $k$, $\{\tau_\alpha < \widetilde{T} < d^+(j_k) < d^+((-\infty, a_\ell])\} \cap H_2^+(j_k) \subset F_2^+(j_k)$, and $\{\tau_\alpha < \widetilde{T} < d^+(j_k) < d^+((-\infty, a_\ell])\} \cap H_3^+(j_k) \subset E_3^+(j_k)$. Therefore, for large $k$,

$$P_{k,\ell}^+ \geq P_z\{\tau_\alpha < \widetilde{T} < d^+(j_k) < d^+((-\infty, a_\ell]), H_2^+(j_k), H_3^+(j_k)\}.$$

We apply the strong Markov property at time $\widetilde{T}$. Since $V(\widetilde{T}) \in I_{a_{\ell+2}} := [a_{\ell+2}, a_{\ell+2} + M]$, we have, for large $k$,

$$P_{k,\ell}^+ \geq P_z\{\tau_\alpha < \widetilde{T} < d^+(j_k), \widetilde{T} < d^+((-\infty, a_\ell]), H_2^+(j_k)\}$$

(6.7)
$$\times \inf_{v \in I_{a_{\ell+2}}} P_v\Big\{d^+(j_k) < d^+((-\infty, a_\ell]),$$

$$\max_{0 \leq x \leq y \leq d^+(j_k)}[V(x) - V(y)] \leq j_k\Big\}.$$

Of course, $\{\tau_\alpha < \widetilde{T}\} = \{\tau_1 < T_1 < \tau_2 < \cdots < T_{\alpha-1} < \tau_\alpha < \widetilde{T}\}$. To estimate $P_z\{\cdots\}$ on the right-hand side, we apply the strong Markov property successively at $\tau_\alpha$, $T_{\alpha-1}$, $\tau_{\alpha-1}, \ldots, T_1$ and $\tau_1$. At time $\tau_\alpha$, we use the following inequality [see (6.4)]: for $v \in [a_{\ell+1} - M, a_{\ell+1})$,

$$P_v\{d^+([a_{\ell+2}, \infty)) < d^+((-\infty, a_\ell])\} \geq \frac{(a_{\ell+1} - M) - a_\ell}{a_{\ell+2} - a_\ell + M} = \frac{2}{7}.$$

At times $\tau_p$ and $T_p$ ($1 \leq p < \alpha$), we use [see (6.4) and (6.5)], respectively, for $v \in [a_{\ell+1} - M, a_{\ell+1})$ and $u \in [a_{\ell+1}, a_{\ell+1} + M]$,

$$P_v\{d^+([a_{\ell+1}, \infty)) < d^+((-\infty, a_\ell])\} \geq \frac{(a_{\ell+1} - M) - a_\ell}{a_{\ell+1} - a_\ell + M} = \frac{1}{2},$$

$$P_u\{d^+((-\infty, a_{\ell+1})) < d^+([a_{\ell+2}, \infty))\} \geq \frac{a_{\ell+2} - (a_{\ell+1} + M)}{a_{\ell+2} - a_{\ell+1} + M} = \frac{1}{2}.$$

Accordingly,

$$P_z\{\tau_\alpha < \widetilde{T} < d^+(j_k), \widetilde{T} < d^+((-\infty, a_\ell]), H_2^+(j_k)\}$$
$$\geq \frac{2/7}{2^{2\alpha-2}} P_z\{\tau_1 < d^+(j_k), H_2^+(j_k)\}.$$

By Donsker's theorem, $\inf_{z \in I_{j_{k-1}}} P_z\{\tau_1 < d^+(j_k), H_2^+(j_k)\}$ is greater than a constant (for large $k$, and uniformly in $\ell$). Thus

$$P_z\{\tau_\alpha < \widetilde{T} < d^+(j_k), \widetilde{T} < d^+((-\infty, a_\ell]), H_2^+(j_k)\} \geq \frac{C_{13}}{2^{2\alpha}} \geq \frac{C_{14}}{k},$$

the last inequality following from the definition of $\alpha := \lfloor \frac{1}{2} \log \log j_k \rfloor$. Plugging this into (6.7) gives that for large $k$,

$$P_{k,\ell}^+ \geq \frac{C_{14}}{k} \inf_{v \in I_{a_{\ell+2}}} P_v\Big\{d^+(j_k) < d^+((-\infty, a_\ell]),$$

$$\max_{0 \leq x \leq y \leq d^+(j_k)}[V(x) - V(y)] \leq j_k\Big\}$$

$$\geq \frac{C_{14}}{k} \inf_{v \in I_{a_{\ell+2}}} P_v\{A_\ell^{(+1)}\} \prod_{2 \leq p \leq 5} \inf_{v \in [((p-4)/2)j_k, ((p-4)/2)j_k + M]} P_v\{A_\ell^{(+p)}\},$$

where
$$A_\ell^{(+1)} := \{d^+([-j_k, \infty)) < d^+((-\infty, a_\ell])\},$$
$$A_\ell^{(+p)} := \left\{d^+\left(\left[\frac{p-3}{2}j_k, \infty\right)\right) < d^+\left(\left(-\infty, \frac{p-5}{2}j_k\right]\right)\right\}, \quad 2 \leq p \leq 5.$$

(The last inequality was obtained by applying the strong Markov property successively at the stopping times $d^+([j_k/2, \infty))$, $d^+([0, \infty))$, $d^+([-j_k/2, \infty))$ and $d^+([-j_k, \infty))$.) It is clear that there exist constants $C_{15} > 0$ and $C_{16} > 0$ such that
$$\inf_{v \in I_{a_{\ell+2}}} P_v\{A_\ell^{(+1)}\} \geq \frac{C_{15}}{j_k},$$
$$\min_{2 \leq p \leq 5} \inf_{v \in [((p-4)/2)j_k, ((p-4)/2)j_k + M]} P_v\{A_\ell^{(+p)}\} \geq C_{16}.$$

Therefore,
$$P_{k,\ell}^+ \geq \frac{C_{14}}{k} \frac{C_{15}}{j_k} (C_{16})^4 =: \frac{C_{17}}{kj_k}.$$

Plugging this into (6.6) gives
$$P_z\{F^+(j_k)\} \geq \left\lfloor \frac{j_k}{3M} \right\rfloor \frac{C_{17}}{kj_k} \geq \frac{C_{18}}{k},$$

which implies (6.3), and completes the proof of (6.1).

6.2. *Proof of* (6.2). We write $V_-(n) := V(-n)$, $\forall n \geq 0$. Let as before $P_z(\cdot) := P(\cdot | V(0) = z)$. Under $P_z$, for $r > z$, we define $d^-(r)$ exactly as in (3.3), that is, $|d^-(r)| := \min\{i \geq 0 : V_-(i) \geq r\}$, whereas for $s < z$, we define
$$|d^-(s)| := \min\{i \geq 0 : V_-(i) \leq s\}.$$

We start with the following estimate: there exist positive constants, denoted by $C_5$ and $C_{19}$, such that

(6.8)
$$\inf_{r \geq 1} P\left\{\sum_{0 \leq x \leq |d^-(r)|} e^{-[V_-(x) - V_-(|b^-(r)|)]} \leq C_5, |d^-(r)| < \left|d^-\left(-\frac{r}{2}\right)\right|\right\}$$
$$\geq C_{19} > 0.$$

This is essentially a consequence of Theorem 2.1 of [2], which is a path decomposition for $(V_-(s), s \leq n)$, when $n$ is deterministic. For more details, we refer to Lemma 3.2 of [9], which, by means of an elementary argument, extends Bertoin's theorem for hitting times. Inequality (6.8) then follows from this lemma via the observation that it is possible to choose $1 + c_{11} > 2c_{13}$ in [9] (notation of [9]) such that when $E_1(t) \cap E_2(r)$ is true (notation of [9]), we have $\min_{0 \le x \le |d^-(r)|} V_-(x) = \min_{0 \le x \le t} V_-(x) \ge -\frac{r}{2}$ (our notation).

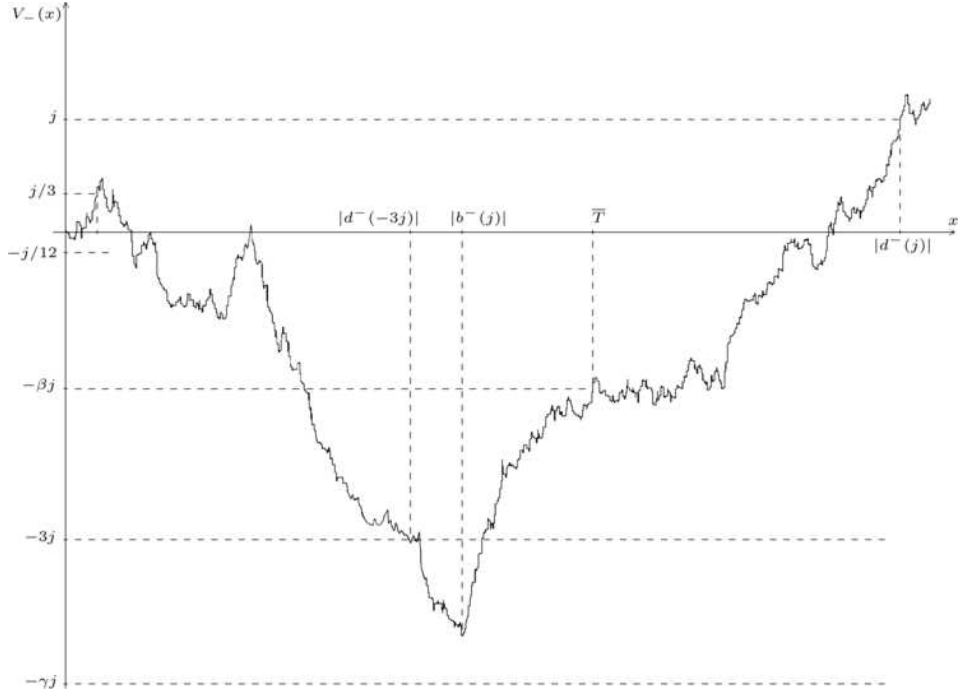

Fig. 2. *Example of $\omega \in \Theta^-(j)$.*

To prove (6.2), we write $\beta := 3 - \frac{1}{1000}$ and $\gamma := 3 + \frac{1}{1000}$, and define

$$\overline{T} := \min\{i \ge |d^-(-3j)| : V_-(i) \ge -\beta j\},$$

$$\widetilde{T} := \min\{i \ge \overline{T} : V_-(i) \le -3j\},$$

$$\Theta^-(j) := \left\{ \left|d^-\left(\frac{j}{3}\right)\right| < \left|d^-\left(-\frac{j}{12}\right)\right| < \overline{T} < |d^-(j)| < \widetilde{T} < |d^-(-\gamma j)| \right\}.$$

See Figure 2 for an example of $\omega \in \Theta^-(j)$.

Recall that $E^-(j) = \bigcap_{i=1}^{5} E_i^-(j)$. Clearly, $E_1^-(j) \cap E_2^-(j) \supset \Theta^-(j)$. Thus

$$E^-(j) \supset \Theta^-(j) \cap E_3^-(j) \cap E_4^-(j) \cap E_5^-(j).$$

Let

$$F_3^-(j) := \left\{ \max_{|d^-(j/3)| \le x \le y \le |d^-(-3j)|} [V_-(y) - V_-(x)] \le \frac{j}{12} \right\},$$

$$F_4^-(j) := \left\{ \frac{j}{3} \le \max_{\overline{T} \le x \le y \le |d^-(j)|} [V_-(x) - V_-(y)] \le j \right\}.$$

Then $E_3^-(j) \supset \Theta^-(j) \cap F_3^-(j)$, and $E_4^-(j) \supset \Theta^-(j) \cap F_4^-(j)$. Thus
$$E^-(j) \supset \Theta^-(j) \cap F_3^-(j) \cap F_4^-(j) \cap E_5^-(j).$$

On $\Theta^-(j) \cap \{|d^-(j)| \le j^3\} \cap \{V_-(|b^-(j)|) \le -3j - j^{1/2}\}$, we have
$$\sum_{x \in [0,|d^-(j)|] \setminus [|d^-(-3j)|,\overline{T}]} e^{-[V_-(x)-V_-(|b^-(j)|)]} \le j^3 e^{-j^{1/2}} \le 1,$$

(for large $j$). Thus $E_5^-(j) \supset F_5^-(j) \cap \Theta^-(j) \cap \{|d^-(j)| \le j^3\} \cap \{V_-(|b^-(j)|) \le -3j - j^{1/2}\}$ (for large $j$), where
$$F_5^-(j) := \left\{ \sum_{|d^-(-3j)| \le x \le \overline{T}} e^{-[V_-(x)-V_-(|\widehat{b}^-(-\beta j)|)]} \le C_5 \right\},$$

and $|\widehat{b}^-(-\beta j)| := \min\{n \ge |d^-(-3j)| : V_-(n) = \min_{x \in [|d^-(-3j)|, \overline{T}]} V_-(x)\}$.

For $j \to \infty$, we have $P\{|d^-(j)| > j^3\} \to 0$ and $P\{V_-(|b^-(j)|) \in (-3j - j^{1/2}, -3j]\} \to 0$. Therefore,

(6.9) $\quad P\{E^-(j)\} \ge P\{\Theta^-(j), F_3^-(j), F_4^-(j), F_5^-(j)\} - o(1),$

where $o(1)$ denotes a term which tends to 0 (when $j \to \infty$). The value of $o(1)$ may vary from line to line.

We apply the strong Markov property at time $\overline{T}$. Since $V_-(\overline{T}) \in I_{-\beta j} := [-\beta j, -\beta j + M]$, this leads to: for large $j$,
$$P\{\Theta^-(j), F_3^-(j), F_4^-(j), F_5^-(j)\} \ge P^{(1)} \inf_{v \in I_{-\beta j}} P_v^{(2)},$$

where
$$P^{(1)} := P\left\{ \left|d^-\left(\frac{j}{3}\right)\right| < \left|d^-\left(-\frac{j}{12}\right)\right| < \overline{T} < |d^-(j)|, \right.$$
$$\left. \overline{T} < |d^-(-\gamma j)|, F_3^-(j), F_5^-(j) \right\},$$
$$P_v^{(2)} := P_v\left\{ |d^-(j)| < |d^-(-3j)|, \frac{j}{3} \le \max_{0 \le x \le y \le |d^-(j)|} [V_-(x) - V_-(y)] \le j \right\}.$$

By Donsker's theorem, $\inf_{v \in I_{-\beta j}} P_v^{(2)} \ge C_{20} > 0$ (for large $j$). Therefore, for large $j$,
$$P\{\Theta^-(j), F_3^-(j), F_4^-(j), F_5^-(j)\} \ge C_{20} P^{(1)}.$$

To obtain a lower bound for $P^{(1)}$, we apply the strong Markov property at time $d^-(-3j)$. Since $V_-(|d^-(-3j)|) \in I_{-3j-M} := [-3j - M, -3j]$, we have
$$P\{\Theta^-(j), F_3^-(j), F_4^-(j), F_5^-(j)\} \ge C_{20} P^{(3)} \inf_{v \in I_{-3j-M}} P_v^{(4)},$$

where
$$P^{(3)} := P\left\{\left|d^-\left(\frac{j}{3}\right)\right| < \left|d^-\left(-\frac{j}{12}\right)\right| < |d^-(-3j)| < |d^-(j)|, F_3^-(j)\right\},$$
$$P_v^{(4)} := P_v\left\{|d^-(-\beta j)| < |d^-(-\gamma j)|, \sum_{x=0}^{|d^-(-\beta j)|} e^{-[V_-(x)-V_-(|b^-(-\beta j)|)]} \leq C_5\right\}.$$

We recall that $|b^-(-\beta j)| := \min\{n \geq 0 : V_-(n) = \min_{x \in [0, |d^-(-\beta j)|]} V_-(x)\}$.

By Donsker's theorem, $P^{(3)}$ is greater than a positive constant (for all large $j$), whereas according to (6.8), $P_v^{(4)} \geq C_{19}$ (for large $j$, uniformly in $v \in I_{-3j-M}$). As a consequence, for large $j$,
$$P\{\Theta^-(j), F_3^-(j), F_4^-(j), F_5^-(j)\} \geq C_{21} > 0.$$
Plugging this into (6.9) completes the proof of (6.2).

**7. A remark.** For any set $A$, let $\xi(n, A) := \sum_{x \in A} \xi(n, x) = \#\{i : 0 \leq i \leq n, X_i \in A\}$.

The recent work of Andreoletti [1] focuses on:
$$Y_n := \inf_{x \in \mathbb{Z}} \min\{k \geq 0 : \xi(n, [x-k, x+k]) \geq an\},$$
where $a \in [0, 1)$ is an arbitrary but fixed constant. In words, $Y_n$ is (half) the minimal size of an interval where the walk hits at least $na$ times in the first $n$ steps.

It is proved in [1] that under (1.1)–(1.3), there exists a constant $c \in (0, \infty)$ such that
$$\liminf_{n \to \infty} Y_n \leq c, \qquad \mathbb{P}\text{-a.s.}$$

A close look at our argument in Section 5 reveals that for some constant $c_* > 0$,
$$(7.1) \qquad \limsup_{n \to \infty} \frac{Y_n}{\log \log \log n} \geq c_*, \qquad \mathbb{P}\text{-a.s.}$$

In fact, the proof of (5.2) shows that, for some constant $C_{22} > 0$, $\max_{x \in [1, d_k^+]} \xi(\tau_k^-, x) \leq C_{22} \frac{\tau_k^-}{\log \log m_k}$ ($\mathbb{P}$-almost surely, for all large $k$; ditto for all the other inequalities stated in this paragraph). In view of (4.6) and (4.16), this implies $\max_{x \in \mathbb{Z}} \xi(\tau_k^-, x) \leq C_{22} \frac{\tau_k^-}{\log \log m_k}$. On the other hand, by (4.15), $\sum_k P_\omega\{\tau_k^- \geq e^{3m_k}, \tau_k^- < \tau(d_k^+)\} < \infty$. Since $\tau_k^- < \tau(d_k^+)$ (Lemma 4.1), we have $\tau_k^- \leq e^{3m_k}$. Thus, $\max_{x \in \mathbb{Z}} \xi(\tau_k^-, x) \leq 2C_{22} \frac{\tau_k^-}{\log \log \log \tau_k^-}$. As a result, (7.1) follows, with $c_* := \frac{1}{2C_{22}}$.

It is, however, not clear whether inequality "$\leq$" would hold in (7.1) with an enlarged value of the constant $c_*$.

**Acknowledgments.** We are grateful to Arvind Singh for help in making the figures in Sections 3 and 6. We wish to thank two anonymous referees for careful reading of the manuscript and for helpful comments.

LABORATOIRE DE PROBABILITÉS
ET MODÈLES ALÉATOIRES
UNIVERSITÉ PARIS VI
4 PLACE JUSSIEU
F-75252 PARIS CEDEX 05
FRANCE
E-MAIL: zhan@proba.jussieu.fr
zindy@ccr.jussieu.fr